\documentclass{birkjour}
 \theoremstyle{definition}
 
 \theoremstyle{remark}

 \numberwithin{equation}{section}
 
 \newcommand{\BR}{{\mathbb R}}

 \newcommand{\BZ}{{\mathbb Z}}

\newcommand{\e}{{\bf e}}

\newcommand{\cl}{C \kern -0.1em \ell}

\begin{document}

%-------------------------------------------------------------------------
% editorial commands: to be inserted by the editorial office
%
%\firstpage{1} \volume{228} \Copyrightyear{2004} \DOI{003-0001}
%\seriesextraline{This is the Concrete Title of this Book\br H.E. R and S.T.C. W, Eds.}
%\firstpage{1}

\title[A note on the discrete CK extension]
 {A note on the discrete Cauchy-Kovalevskaya extension}

%----------Author 1
%\author[Birkh\"auser]{Birkh\"{a}user Publishing Ltd.}

%\address{%
%Viaduktstr. 42\\
%P.O. Box 133\\
%CH 4010 Basel\\
%Switzerland}

%\email{info@birkhauser.ch}

%\thanks{This work was completed with the support of our
%----------Author 2
\author{Nelson Faustino}
\address{Universidade Federal do ABC \br
	Avenida dos Estados, 5001 \br
	CEP 09210-580, Santo Andr\'e \br
	Brazil}
\email{nelson.faustino@ufabc.edu.br}



\date{\today}

\begin{abstract}
In this paper we exploit the umbral calculus framework to reformulate the so-called discrete Cauchy-Kovalevskaya extension in the scope of hypercomplex variables. 
The key idea is to consider not only formal power series representation for the underlying solution, but also integral representations for the Chebyshev polynomials of first and second kind by means of its Cauchy principal values.
It turns out that the resulting integral representation associated to our toy problem is a space-time Fourier type inversion formula. Moreover, with the aid of some Laplace transform identities involving the generalized Mittag-Leffler function we are able to establish a link with a Cauchy problem of differential-difference type.
\end{abstract}

%%% ----------------------------------------------------------------------
\maketitle
%%% ----------------------------------------------------------------------
%\tableofcontents

\section{Introduction}

Apart the big advances of Augustin-Louis Cauchy
(cf.~\cite[Subsection 2.3]{Cooke84}) and Karl
Weierstra\ss~(cf.~\cite[Subsection 2.4]{Cooke84}), Sonya
Kovalevskaya's approach (cf.~\cite[Subsection 2.5]{Cooke84}) may be considered
as one of the major landmarks obtained on the {\it nineteenth-century} for
the study of analytic representations for the solutions of partial differential equations (PDE's).

Its uselfulness in the field of hypercomplex variables is undisputed, largely due to the pioneering work of Fransiscus Sommen 
\cite{Sommen81}. An immediate consequence of such kind of construction is that variants of Kovalevskaya's results -- namely the
ones associated to PDE's in the 'normal form'
(cf.~\cite[p.~32]{Cooke84})-- may be used successively as a tool
to generate null solutions for Dirac-type operators, the so-called
monogenic functions. We refer to \cite[pp.~265-268]{DSS92} for the
generalized theory in the context of hypercomplex variables, and to
\cite[Chapter III]{DSS92} for a wide class of examples involving spherical
monogenics and other generating classes of special functions. For a fractional calculus counterpart of Kovalevskaya's result involving Gegenbauer polynomials we refer to \cite{Vieira17}.

The discrete counterpart for the Cauchy-Kovalevskaya's approach was already investigated by De Ridder et al (cf.~\cite{RSS10}) and Constales \& De Ridder (cf.~\cite{CR14}) using slightly different approaches.
In the first approach (2010), the authors mimic the continuous counterpart by means of Taylor series expansions in interplay with the so-called \textit{skew-Weyl relations}, formely introduced on the paper \cite{RSKS10}. In the second approach (2014) the authors provided an alternative construction for \cite{RSS10}, using solely operational identities associated to the Chebyshev polynomials of first and second kind (cf.~\cite[p. 170]{MasonChebyshev93}):
\begin{eqnarray}
\label{ChebyshevPolynomials}T_k(\lambda)=\cos(k\cos^{-1}(\lambda)) & \mbox{resp.}& \displaystyle U_{k-1}(\lambda)=\dfrac{\sin(k\cos^{-1}(\lambda))}{\sqrt{1-\lambda^2}}.
\end{eqnarray}

In the present paper we reformulate the approaches considered on both papers. The construction stated in section \ref{CKsection} makes use of Roman's umbral calculus approach \cite{Roman84} and of the discrete Fourier theory considered by G\"urlebeck and Spr\"o\ss ig on their former book \cite{GuerlebeckSproessig97},
following up author's recent contribution \cite{FaustinoRelativistic18}. 
In section \ref{PlaneWaveSection} it was obtained an explicit \textit{space-time} integral representation for the Cauchy-Kovalevskaya extension on the momentum space, based on Cauchy principal value representations for (\ref{ChebyshevPolynomials}), and on a Laplace transform identity associated to the \textit{generalized Mittag-Leffler function} (cf.~\cite[p.~21]{SamkoEtAl93})

\begin{eqnarray} \label{MittagLeffler}E_{\alpha,\beta}(z)=\sum_{k=0}^{\infty}\frac{z^k}{\Gamma(\beta+\alpha k)}, &\mbox{for}& \Re(\alpha)>0~\&~\Re(\beta)>0.
\end{eqnarray}

With such construction we provide a further effort toward opening a wider scenario on the theory of discrete hypercomplex variables in interplay with fractional integro-differential type operators (cf.~\cite[Chapter 5]{SamkoEtAl93}), with the hope that this may contribute to further developments on the field towards the discrete analogues for the \textit{Littlewood-Paley-Stein theory} (cf. \cite{SteinW2000,MSteinW02,Pierce09}). 

\section{Problem setup}

The approach proposed throughout this paper is centered around the solution of a first order difference-difference Cauchy problem
on the space-time lattice $$h\BZ^n \times \tau \BZ_{\geq 0}:=\left\{(x,t)\in \BR^n \times [0,\infty)~:~\frac{x}{h}\in \BZ^n~~\&~~\frac{t}{\tau}\in \BZ \right\}$$ involving the finite difference discretization $D_h$ for the Dirac operator over the Clifford algebra $\cl_{n,n}$, considered by the author in \cite{FaustinoKGordonDirac16,FaustinoMMAS17}.

To motivate our approach let us take at a first glance a close look to the following second order Cauchy problem on $h\BZ^n \times \tau \BZ_{\geq 0}$: 
\begin{eqnarray}
\label{KleinGordonh} \left\{\begin{array}{lll} 
\dfrac{\Psi(x,t+\tau)+\Psi(x,t-\tau)-2\Psi(x,t)}{\tau^2}= \Delta_h \Psi(x,t) & \mbox{for} & (x,t)\in
h\BZ^n \times \tau \BZ_{\geq 0}
\\ \ \\
\Psi(x,0)=\Phi_0(x) & \mbox{for} & x\in h\BZ^n\\ \ \\
\left[L_t \Psi(x,t)\right]_{t=0}=\Phi_1(x) & \mbox{for} & x\in h\BZ^n
\end{array}\right.
\end{eqnarray}

The above problem corresponds to a difference-difference discretization of the Cauchy problem of Klein-Gordon type the \textit{massless limit} $m\rightarrow 0$ (cf.~\cite[section 4.]{Rabin82}).

Here
\begin{eqnarray}
\label{discreteLaplacian}
\displaystyle \Delta_h \Psi(x,t)=\sum_{j=1}^n
\frac{\Psi(x+h\e_j,t)+\Psi(x-h\e_j,t)-2\Psi(x,t)}{h^2}.
\end{eqnarray}
corresponds to the finite difference action of the \textit{discrete Laplacian} on the space-time lattice $h\BZ^n \times \tau \BZ_{\geq 0}$ (cf. \cite[subsection 1.2.]{FaustinoRelativistic18}), whereas $L_t$ is a finite difference operator satisfying the second-order constraint
$$
L_t(L_t\Psi(x,t))=\dfrac{\Psi(x,t+\tau)+\Psi(x,t-\tau)-2\Psi(x,t)}{\tau^2}.
$$ 

For the case where $L_t$ is a finite difference operator defined as $$L_t\Psi(x,t)=\dfrac{\Psi\left(x,t+\frac{\tau}{2}\right)-\Psi\left(x,t-\frac{\tau}{2}\right)}{\tau},$$ it was depicted by the author in his recent contribution (cf.~\cite[subsection 4.1.]{FaustinoRelativistic18}) that the solution of the difference-difference evolution problem (\ref{KleinGordonh}) may be derived from the Cauchy principal values representations for the Chebyshev polynomials of first and second kind, $T_k(\lambda)$ resp. $U_{k-1}(\lambda)$
(cf.~ \cite[subsection 4.1,p.~173]{MasonChebyshev93}):
\begin{eqnarray}
\label{ChebyshevIntegrals}
\begin{array}{lll}
\displaystyle T_{k}(\lambda)&=&\displaystyle -\frac{1}{\pi}\int_{-1}^{1} \left(1-s^2\right)^{\frac{1}{2}}\frac{U_{k-1}(s)}{s-\lambda} ds\\ \ \\
\displaystyle U_{k-1}(\lambda)&=&\displaystyle\frac{1}{\pi}\int_{-1}^{1} \left(1-s^2\right)^{-\frac{1}{2}}\frac{T_{k}(s)}{s-\lambda} ds.
\end{array}
\end{eqnarray}

As studied in depth by the author in \cite{FaustinoMonomiality14}, the above construction yields as a natural exploitation of an abstract framework involving manipulations of \textit{shift-invariant operators} that admits formal series expansions in terms of the time derivative $\partial_t$: 
$$ 	L_t=\sum_{k=1}^\infty b_k \dfrac{\left(\partial_t\right)^k}{k!},~~~\mbox{with}~~~b_k=[(L_t)^k t^k]_{t=0}.$$

Such construction looks similar to \cite{CR14} (cf.~\cite[Remark 4.1]{FaustinoRelativistic18}), except that the $-iL_t$ operator ($i=\sqrt{-1}$) -- appearing quite often on Dirac's equation (cf.~\cite[subsection 5.]{Rabin82}) -- is replaced by finite difference operators of the form
\begin{eqnarray}
\label{FiniteDiffCR14}\nabla_\tau \Psi(x,t)=\e_- \frac{\Psi(x,t)-\Psi(x,t-\tau)}{\tau}-\e_+ \frac{\Psi(x,t+\tau)-\Psi(x,t)}{\tau},
\end{eqnarray}
for two given Witt basis on the time direction, $\e_+$ resp. $\e_-$, satisfying the set of constraints
\begin{eqnarray*}
	\left(\e_+\right)^2=\left(\e_-\right)^2=0 & \& & \e_+ \e_-+\e_-\e_+=1.
\end{eqnarray*}

The Cauchy problem that we propose to study here is of the form
\begin{eqnarray}
\label{DiracCK} \left\{\begin{array}{lll} 
\e_0 \dfrac{\Psi(x,t+\tau)-\Psi(x,t-\tau)}{2\tau}+\e_{2n+1}\dfrac{2\Psi(x,t)-\Psi(x,t+\tau)-\Psi(x,t-\tau)}{2\tau} \\ \ \\
=- D_h \Psi(x,t),~~~~\mbox{for}~~~(x,t)\in h\BZ^n \times \tau \BZ_{\geq 0} \\ \ \\
\Psi(x,0)=\Phi_0(x),~~\mbox{for}~~~x\in h\BZ^n
\end{array}\right..  
\end{eqnarray}

Here and elsewhere
\begin{eqnarray}
\label{DiracEqh}
\begin{array}{lll}
D_h \Psi(x,t)&=&\displaystyle \sum_{j=1}^n\e_j\frac{\Psi(x+h \e_j,t)-\Psi(x-h\e_j,t)}{2h}+ \\
&+&	\displaystyle \sum_{j=1}^n\e_{n+j}\frac{2\Psi(x,t)-\Psi(x+h \e_j,t)-\Psi(x-h \e_j,t)}{2h}
\end{array}
\end{eqnarray}
stands for the finite difference discretization for the Dirac operator, whereas $\e_0$ and $\e_{2n+1}$ are two Clifford generators satisfying 
\begin{eqnarray*}
	(\e_0)^2=-1 & \& & (\e_{2n+1})^2=+1.
\end{eqnarray*}

The resulting Cauchy problem formulation over the space-time lattice $h\BZ^n \times \tau \BZ_{\geq 0}$ mix first order approximations for the space-time derivatives, $\partial_{t}$ and $\partial_{x_j}$ respectively, on the $\e_j-$direction ($j=0,1,\ldots,n$) plus $n+1$ second order perturbation terms, satisfying the asymptotic conditions
\begin{eqnarray*}
	\dfrac{2\Psi(x,t)-\Psi(x,t+\tau)-\Psi(x,t-\tau)}{2\tau} &=&-\frac{\tau}{2}\partial_{t}^2\Psi(x,t)+O(\tau^3)\\ \ \\	
	\frac{2\Psi(x,t)-\Psi(x+h \e_j,t)-\Psi(x-h \e_j,t)}{2h}&=&-\frac{h}{2}\partial_{x_j}^2\Psi(x,t)+O(h^3).
\end{eqnarray*}

Throughout this paper we shall assume that $\e_0$ and $\e_{2n+1}$ together with $\e_1,\e_2,\ldots,\e_n$, $\e_{n+1},\e_{n+2}\,\ldots,\e_{2n}$ are the generators of the Clifford algebra $\cl_{n+1,n+1}$, satisfying 
\begin{eqnarray}
\label{CliffordBasis}
\begin{array}{lll}
\e_j \e_k+ \e_k \e_j=-2\delta_{jk}, & 0\leq j,k\leq n \\
\e_{j} \e_{n+k}+ \e_{n+k} \e_{j}=0, & 0\leq j\leq n ~~\&~~  1\leq k\leq n+1\\
\e_{n+j} \e_{n+k}+ \e_{n+k} \e_{n+j}=2\delta_{jk}, & 1\leq j,k\leq
n+1.
\end{array}
\end{eqnarray}

Under the canonical isomorphism $\cl_{n+1,n+1}\cong \mbox{End}(\cl_{0,n+1})$ (see e.g. \cite[Chapter 4]{VazRoldao16} for further details), the equivalence between our formulation and the formulation considered by Constales and De Ridder in \cite{CR14} is rather obvious (cf.~\cite[subsection 2.3]{FaustinoKGordonDirac16}).

In particular, under the canonical identifications
\begin{eqnarray*}
	\e_0 \leftrightarrow \e_--\e_+& \& & \e_{2n+1}\leftrightarrow \e_-+\e_+
\end{eqnarray*}
we find that the finite difference operators
$$\e_0 \dfrac{\Psi(x,t+\tau)-\Psi(x,t-\tau)}{2\tau}+\e_{2n+1}\dfrac{2\Psi(x,t)-\Psi(x,t+\tau)-\Psi(x,t-\tau)}{2\tau}$$
and (\ref{FiniteDiffCR14}) are indeed equivalent.

\section{The discrete Cauchy-Kovalevskaya approach explained}\label{CKsection}

\subsection{Starting with a formal power series expansion}\label{EGFsection}

The formal construction of the discrete Cauchy-Kovalevskaya extension proceeds as follows: 

Let us take first a close look for the equation
$$		\e_0 \dfrac{\Psi(x,t+\tau)-\Psi(x,t-\tau)}{2\tau}+\e_{2n+1}\dfrac{2\Psi(x,t)-\Psi(x,t+\tau)-\Psi(x,t-\tau)}{2\tau}
=- D_h \Psi(x,t).$$ 

We notice here that the factorization property $\left(D_h\right)^2=-\Delta_h$ (cf.~\cite[Proposition 2.1]{FaustinoKGordonDirac16}) together with the identities $(\e_0)^2=-1$ and $(\e_{2n+1})^2=+1$ lead to 
\begin{eqnarray}
\label{discreteHarmonic}
\begin{array}{lll}
-\Delta_h \Psi(x,t)&=&-D_h\left(-D_h\Psi(x,t)\right)\\ \ \\
&=&\dfrac{2\Psi(x,t)-\Psi(x,t+\tau)-\Psi(x,t-\tau)}{\tau^2}.
\end{array}
\end{eqnarray}

From the combination of the above two identities, we end up with
$$	\e_0 \dfrac{\Psi(x,t+\tau)-\Psi(x,t-\tau)}{2\tau}-\e_{2n+1}\frac{\tau}{2}\Delta_h \Psi(x,t)
=- D_h \Psi(x,t),
$$

Moreover, from the set of identities $\left(\e_0\right)^2=-1$ and $\e_0\e_{2n+1}=-\e_{2n+1}\e_0$ involving the Clifford generators $\e_0$ and $\e_{2n+1}$ of $\cl_{n+1,n+1}$, there holds
\begin{eqnarray}
\label{DiracEqhtau}
\Psi(x,t+\tau)-\Psi(x,t-\tau)=\left(2\tau\e_0 D_h+\tau^2\e_{2n+1}\e_0\Delta_h\right) \Psi(x,t).
\end{eqnarray}

Here we recall that the left-hand side of (\ref{DiracEqhtau}) admits the formal series expansion (cf.~\cite[Example 2.3.]{FaustinoMonomiality14})
$$
\Psi(x,t+\tau)-\Psi(x,t-\tau)=2\sinh(\tau\partial_t)\Psi(x,t).
$$

Similarly to \cite[subsection 3.2.]{FaustinoRelativistic18}, one can conclude from the identity (\ref{DiracEqhtau}) that the solution of the Cauchy problem (\ref{DiracCK}) may be represented through the {\it Exponential Generating Function} (EGF, for short) type expansion
\begin{eqnarray}
\label{CKseriesExpansion}
\Psi(x,t)=\sum_{k=0}^{\infty} \frac{G_k(t;-\tau,2\tau)}{k!}\left(2\tau\e_0 D_h+\tau^2\e_{2n+1}\e_0\Delta_h\right)^k \Phi_0(x).
\end{eqnarray}

On the above formula $G_k(t;-\tau,2\tau)$ ($k \in \mathbb{N}_0$) denote the Gould polynomials -- the Sheffer sequence associated to the \textit{delta operator} $2\sinh(\tau \partial_t)$ (cf.~\cite[subsection 1.4.]{Roman84}).

We recall here that (\ref{CKseriesExpansion}) slightly differs from the one proposed in \cite[pp.~1469 \& 1470]{RSS10}. The main idea here to rid of the need of imposing skew-Weyl constraints (cf.~\cite[section 3]{RSKS10}) was the link between the solution of the discrete Cauchy problem (\ref{DiracCK}) with the solutions of the {discrete harmonic type equation} (\ref{discreteHarmonic}) that lead us to a first order time-evolution equation (\ref{DiracEqh}) encoded by the {\it delta operator} $2\sinh(\tau \partial_t)$. 

\subsection{Operational Identities}\label{OperationalSection}

In the previous section we have used the EGF endowed by the Sheffer sequence associated to the \textit{delta operator} $2\sinh(\tau\partial_t)$, to obtain a formal series representation for the solution of the Cauchy problem (\ref{DiracEqh}).

Now we will see how the factorization of the finite difference operator $2\tau\e_0 D_h+\tau^2\e_{2n+1}\e_0\Delta_h$ together with rather simple properties involving umbral calculus techniques lead to some interesting operational identities.

We start to observe that the set of relations (\ref{CliffordBasis}) associated to the Clifford generators of $\cl_{n+1,n+1}$ lead to the set of anti-commuting relations 
\begin{eqnarray*}
	\e_0 D_h+D_h \e_0=0 \\
	\e_0 D_h (\e_{2n+1}\e_0 \Delta_h)+(\e_{2n+1}\e_0 \Delta_h)\e_0 D_h=0,
\end{eqnarray*}
and hence, to the factorization property
\begin{eqnarray*}
	\label{factorizationDhtau}\left(2\tau\e_0 D_h+\tau^2\e_{2n+1}\e_0\Delta_h\right)^2=-4\tau^2\Delta_h+ \tau ^4 \Delta_h^2.
\end{eqnarray*}

Thereby, the iterated powers $\left(2\tau\e_0 D_h+\tau^2\e_{2n+1}\e_0\Delta_h\right)^{k}$ appearing on the right-hand side of the infinite summand (\ref{CKseriesExpansion}) may be splitted on even ($k=2m$) and odd powers ($k=2m+1$) as follows:
\begin{eqnarray}
\label{powersDhtauk}
\small{\begin{array}{lll}
\left(2\tau\e_0 D_h+\tau^2\e_{2n+1}\e_0\Delta_h\right)^{2m}&=&\left(-4\tau^2\Delta_h+ \tau ^4 \Delta_h^2\right)^m
\\ 
\left(2\tau\e_0 D_h+\tau^2\e_{2n+1}\e_0\Delta_h\right)^{2m+1}&=&\left(-4\tau^2\Delta_h+ \tau ^4 \Delta_h^2\right)^m\left(2\tau\e_0 D_h+\tau^2\e_{2n+1}\e_0\Delta_h\right).
\end{array}}
\end{eqnarray}

By the substitution of (\ref{powersDhtauk}) into (\ref{CKseriesExpansion}) we recognize, after some umbral calculus manipulations involving the inverse $L^{-1}(s)=\frac{1}{\tau}\sinh^{-1}\left(\frac{s}{2}\right)$ of the function $L(s)=2\sinh(\tau s)$ (cf.~\cite[Appendix A]{FaustinoRelativistic18}) that
\begin{eqnarray}
\label{CKseriesExpansionUmbral}
\begin{array}{lll}
\Psi(x,t)&=&
\displaystyle\sum_{m=0}^{\infty} \frac{2^{2m}G_{2m}(t;-\tau,2\tau)}{(2m)!}\left(-\tau^2\Delta_h+ \frac{\tau^4}{4} \Delta_h^2\right)^m\Phi_0(x)+\\ \ \\ &+&\displaystyle \sum_{m=0}^{\infty} \frac{2^{2m+1}G_{2m+1}(t;-\tau,2\tau)}{(2m+1)!}\left(-\tau^2\Delta_h+ \frac{\tau^4}{4} \Delta_h^2\right)^m \times \\ \ \\ &\times & \left[\tau\e_0 D_h\Phi_0(x)+\dfrac{\tau^2}{2}\e_{2n+1}\e_0\Delta_h\Phi_0(x)\right] \\ \ \\
&=&\cosh\left(\dfrac{t}{\tau}\sinh^{-1}\left(\sqrt{-\tau^2\Delta_h+ \dfrac{\tau^4}{4} \Delta_h^2}\right)\right)\Phi_0(x)+\\ \ \\
&+&\dfrac{\sinh\left(\dfrac{t}{\tau}\sinh^{-1}\left(\sqrt{-\tau^2\Delta_h+ \dfrac{\tau^4}{4} \Delta_h^2}\right)\right)}{\sqrt{-\tau^2\Delta_h+ \dfrac{\tau^4}{4} \Delta_h^2}}\times \\ \ \\ &\times &\left[\tau\e_0 D_h\Phi_0(x)+\dfrac{\tau^2}{2}\e_{2n+1}\e_0\Delta_h\Phi_0(x)\right].
\end{array}
\end{eqnarray}

\subsection{The discrete Fourier analysis approach}

In analogy to \cite[section 3.]{FaustinoRelativistic18}, the computation of the solution for the time-evolution problems may be obtained with the aid 
of the discrete Fourier transform (cf.~\cite[subsection 5.2.1]{GuerlebeckSproessig97})
\begin{eqnarray}
\label{discreteFh}
(\mathcal{F}_{h} \Psi)(\xi,t)&=&\left\{\begin{array}{lll}
\displaystyle \frac{h^n}{\left(2\pi\right)^{\frac{n}{2}}}\displaystyle 
\sum_{x\in h\BZ^n}\Psi(x,t)e^{i x \cdot \xi} & \mbox{for} & \xi\in Q_h
\\ \ \\
0 & \mbox{for} & \xi\in \BR^n \setminus \left(-\frac{\pi}{h},\frac{\pi}{h} \right]^n,
\end{array}\right.
\end{eqnarray}
where $Q_h=\left(-\frac{\pi}{h},\frac{\pi}{h} \right]^n$ stands for the \textit{Brioullin representation} of the $n-$torus $\BR^n/\frac{2\pi}{h}\BZ^n$ (cf.~\cite[p.~
324]{Rabin82}).
Here we recall that the map
$\mathcal{F}_h:\ell_2(h\BZ^n)\rightarrow L_2(Q_h)$ is an isometry so that
\begin{eqnarray}
\label{FourierTransform} (\mathcal{F}^{-1}_h
\Phi)(x,t)=\frac{1}{\left(2\pi\right)^{\frac{n}{2}}}\int_{Q_h}
\Phi(\xi,t)e^{-ix\cdot \xi }d\xi.
\end{eqnarray}

Using the fact that the (Clifford) constants 
\begin{eqnarray}
\label{FourierMultipliers}
\begin{array}{lll}
d_h(\xi)^2&=&\displaystyle \sum_{j=1}^{n}\frac{4}{h^2}\sin^2\left(\frac{h\xi_j}{2}\right) \\ \ \\ 
\textbf{z}_{h}(\xi)
&=&\displaystyle \sum_{j=1}^n -i\e_j\dfrac{\sin(h\xi_j)}{h}+\sum_{j=1}^n \e_{n+j}\dfrac{1-\cos(h\xi_j)}{h}
\end{array}
\end{eqnarray} 
are the Fourier multipliers of $\mathcal{F}_h\circ (-\Delta_h)\circ\mathcal{F}_h^{-1}$ resp. $\mathcal{F}_h\circ D_h\circ \mathcal{F}_h^{-1}$, we get that the EGF operational identity (\ref{CKseriesExpansionUmbral}) on the momentum space $Q_h \times \tau \BZ_{\geq 0}$ reads as 
\begin{eqnarray}
\label{CKseriesFourier}
\begin{array}{lll}
\mathcal{F}_h\Psi(\xi,t)&=&\cos\left(\dfrac{t}{\tau}\sin^{-1}\left(\sqrt{\tau^2d_h(\xi)^2+ \dfrac{\tau^4}{4} d_h(\xi)^4}\right)\right)\mathcal{F}_h\Phi_0(\xi)+\\ \ \\
&+&\dfrac{\sin\left(\dfrac{t}{\tau}\sin^{-1}\left(\sqrt{\tau^2d_h(\xi)^2+ \dfrac{\tau^4}{4} d_h(\xi)^4}\right)\right)}{\sqrt{\tau^2d_h(\xi)^2+ \dfrac{\tau^4}{4} d_h(\xi)^4}}\times \\ &\times&\left[\tau\e_0 {\bf z}_h(\xi)-\dfrac{\tau^2}{2}\e_{2n+1}\e_0d_h(\xi)^2\right]\mathcal{F}_h\Phi_0(\xi),
\end{array}
\end{eqnarray}

On the other hand, with the aid of the fundamental trigonometric identity 
\begin{eqnarray*}
	\sin^{-1}(s)=\cos^{-1}\left(\sqrt{1-s^2}\right)&\mbox{for}&0\leq s\leq 1,
\end{eqnarray*}
we realize that the condition $\displaystyle d_h(\xi)^2 \leq \frac{2 }{\tau^2}(\sqrt{2}-1)$ that yields from the constraint
\begin{eqnarray*}
	0\leq \tau^2 d_h(\xi)^2+\dfrac{\tau^4}{4} d_h(\xi)^4  \leq 1 
\end{eqnarray*}
allows us to represent (\ref{CKseriesFourier}) in terms of the Chebyshev polynomials of first and second kind, $T_k(\lambda)$ resp. $U_{k-1}(\lambda)$, defined viz (\ref{ChebyshevPolynomials}).
Namely, one has
\begin{eqnarray}
\label{CKseriesFourier}
\begin{array}{lll}
\mathcal{F}_h\Psi(\xi,t)&=&T_{\frac{t}{\tau}}\left(\sqrt{1-\tau^2d_h(\xi)^2- \dfrac{\tau^4}{4} d_h(\xi)^4}\right)\mathcal{F}_h\Phi_0(\xi)~+\\ \ \\
&+&
U_{\frac{t}{\tau}-1}\left(\sqrt{1-\tau^2d_h(\xi)^2- \dfrac{\tau^4}{4} d_h(\xi)^4}\right)\times \\ \ \\
&\times & \left[\tau\e_0 {\bf z}_h(\xi)-\dfrac{\tau^2}{2}\e_{2n+1}\e_0d_h(\xi)^2\right]\mathcal{F}_h\Phi_0(\xi).
\end{array}
\end{eqnarray}

Finally, in the view of the discrete Fourier inversion properties (cf.~\cite[p.~247]{GuerlebeckSproessig97}) we get
that the formal solution (\ref{CKseriesExpansion}) is given by the \textit{discrete convolution formula}
\begin{eqnarray}
\label{discreteConvolutionformula}
\Psi(x,t)&=&\sum_{y\in h\BZ^n} h^n {\bf K}_\tau(x-y,t)\Phi_0(y), 
\end{eqnarray}
involving the kernel function
\begin{eqnarray}
\label{ChebyshevKernel}
\begin{array}{lll}
{\bf K}_\tau(x,t)&=& \dfrac{1}{(2\pi)^{\frac{n}{2}}}\displaystyle \int_{Q_h}T_{\frac{t}{\tau}}\left(\sqrt{1-\tau^2d_h(\xi)^2~- \dfrac{\tau^4}{4} d_h(\xi)^4}\right) e^{-i x\cdot \xi}d\xi~+\\ \ \\
&+& \dfrac{1}{(2\pi)^{\frac{n}{2}}}\displaystyle \int_{Q_h}U_{\frac{t}{\tau}-1}\left(\sqrt{1-\tau^2d_h(\xi)^2- \dfrac{\tau^4}{4} d_h(\xi)^2}\right) \times \\ 
&\times & \left[\tau\e_0 {\bf z}_h(\xi)-\dfrac{\tau^2}{2}\e_{2n+1}\e_0d_h(\xi)^2\right] e^{-i x\cdot \xi}d\xi. 
\end{array}
\end{eqnarray}

\section{Integral Representation formulae}\label{PlaneWaveSection}

\subsection{A Fourier representation formula for $K_\tau(x,t)$ on the space-time}\label{ChebyshevSection}

We start this subsection by exploiting a
quite unexpected hypersingular integral representation for the kernel function (\ref{ChebyshevKernel}) by the Cauchy principal value representations (\ref{ChebyshevIntegrals}) for the Chebyshev polynomials of the first and second kind.
It turns out that it does not require any substantially new ideas: We reformulate only the construction considered in \cite[subsection 5.1]{FaustinoRelativistic18}.

On the sequel we will make use of the factorization identities (cf.~\cite[subsection 2.1.]{FaustinoRelativistic18})
\begin{center}
	${\bf z}_h(\xi)^2=d_h(\xi)^2$ \&  ${\bf z}_h(\xi)^4=d_h(\xi)^4$
\end{center}
involving the Fourier multipliers (\ref{FourierMultipliers}) of $\mathcal{F}_h\circ (-\Delta_h)\circ\mathcal{F}_h^{-1}$ resp. $\mathcal{F}_h\circ D_h\circ \mathcal{F}_h^{-1}$.

First, we observe that for the change of variable $s=\cos(\omega \tau)$, with $0\leq \omega\leq \dfrac{\pi}{\tau}$, 
it easily follows by straightforward integral simplifications involving parity arguments that the Cauchy principal value representations (\ref{ChebyshevIntegrals}) for the Chebyshev polynomials (\ref{ChebyshevPolynomials}) may be rewritten as (cf.~\cite[p.~173]{MasonChebyshev93})
\begin{eqnarray}
\label{ChebyshevIntegralsTrigonometric}
\begin{array}{lll}
T_{\frac{t}{\tau}}(\lambda)=&\displaystyle- \dfrac{\tau}{\pi} \int_{0}^{\frac{\pi}{\tau}} \frac{\sin(\omega \tau)}{\cos(\omega \tau)-\lambda} ~\sin(\omega t)d\omega =&\displaystyle \dfrac{\tau}{2\pi} \int_{-\frac{\pi}{\tau}}^{\frac{\pi}{\tau}} \frac{-i\sin(\omega \tau)}{\cos(\omega \tau)-\lambda}~e^{-i\omega t} d\omega \\ \ \\
U_{\frac{t}{\tau}-1}(\lambda)=&\displaystyle \dfrac{\tau}{\pi} \int_{0}^{\frac{\pi}{\tau}} \frac{1}{\cos(\omega \tau)-\lambda}~\cos(\omega t) d\omega =&\displaystyle \dfrac{\tau}{2\pi} \int_{-\frac{\pi}{\tau}}^{\frac{\pi}{\tau}} \frac{1}{\cos(\omega \tau)-\lambda} e^{-i\omega t}d\omega.	
\end{array}
\end{eqnarray}

Hence, by inserting the integral identities (\ref{ChebyshevIntegralsTrigonometric}) on (\ref{ChebyshevKernel}), there follows through the substitution $\lambda=\sqrt{1-\tau^2{\bf z}_h(\xi)^2- \frac{\tau^4}{4} {\bf z}_h(\xi)^4}$ that
\begin{eqnarray*}
	\label{KtauChebyshev}
{\bf K}_\tau(x,t)= \nonumber \\
= \frac{\tau}{(2\pi)^{\frac{n}{2}+1}}\int_{Q_h} \int_{-\frac{\pi}{\tau}}^{\frac{\pi}{\tau}}\frac{-i\sin(\omega \tau)+\tau\e_0 {\bf z}_h(\xi)-\frac{\tau^2}{2}\e_{2n+1}\e_0{\bf z}_h(\xi)^2}{\cos(\omega\tau)-\sqrt{1-\tau^2{\bf z}_h(\xi)^2- \frac{\tau^4}{4} {\bf z}_h(\xi)^4}}~e^{-i(\omega t+x
		\cdot \xi)}d\omega d\xi.
\end{eqnarray*}

The above formula provides us an integral representation of (\ref{ChebyshevKernel}) over $Q_h \times \left(-\frac{\pi}{\tau},\frac{\pi}{\tau}\right]$. This representation, that fulfils for values of $d_h(\xi)^2$ satisfying the condition $$\displaystyle d_h(\xi)^2 \leq \frac{2 }{\tau^2}(\sqrt{2}-1)$$ is nothing else than a \textit{space-time Fourier inversion type formula} encoded by the
solution of (\ref{DiracCK}) on the momentum space $Q_h \times \left(-\frac{\pi}{\tau},\frac{\pi}{\tau}\right]$.

\subsection{Towards a fractional integro-differential type representation}

We end this note by providing an alternative way to describe (\ref{ChebyshevKernel}) by means of a fractional integral representation of Bessel type (cf.~\cite[part 27 of Chapter 5]{SamkoEtAl93}) for the term 
\begin{eqnarray}
\label{IntegrandChebyshev}
\frac{-i\sin(\omega \tau)+\tau\e_0 {\bf z}_h(\xi)-\frac{\tau^2}{2}\e_{2n+1}\e_0{\bf z}_h(\xi)^2}{\cos(\omega\tau)-\sqrt{1-\tau^2{\bf z}_h(\xi)^2- \frac{\tau^4}{4} {\bf z}_h(\xi)^4}}.
\end{eqnarray}

To do so, we make use of the Laplace transform identity (cf.~\cite[p.~21]{SamkoEtAl93}) 
\begin{eqnarray}
\label{LaplaceIdentityMittagLeffler}
\int_{0}^\infty e^{p\lambda^2} p^{\beta-1}E_{\alpha,\beta}\left(s p^{\alpha}\right)dp= \dfrac{\lambda^{-2\beta}}{1-s \lambda^{-2\alpha}},~~\mbox{for}~~ \Re(\lambda^2)>|s|^{\frac{1}{\alpha}}~\&~\Re(\beta)>0
\end{eqnarray}
involving the \textit{generalized Mittag-Leffler function} $E_{\alpha,\beta}$ defined viz equation (\ref{MittagLeffler}), to obtain a regularization for the term (\ref{IntegrandChebyshev}).

We recall here that for the substitutions
\begin{eqnarray*}
	s=\cos(\omega\tau),~~~~ 
	\lambda=\sqrt{1-\tau^2{\bf z}_h(\xi)^2- \frac{\tau^4}{4} {\bf z}_h(\xi)^4} & \& & \alpha=\beta=\frac{1}{2} 
\end{eqnarray*}
on (\ref{LaplaceIdentityMittagLeffler}) it readily follows that (\ref{IntegrandChebyshev}) admits the following integral representation
\begin{eqnarray*}
	\frac{-i\sin(\omega \tau)+\tau\e_0 {\bf z}_h(\xi)-\frac{\tau^2}{2}\e_{2n+1}\e_0{\bf z}_h(\xi)^2}{\cos(\omega\tau)-\sqrt{1-\tau^2{\bf z}_h(\xi)^2- \frac{\tau^4}{4} {\bf z}_h(\xi)^4}}= \\ \ \\
	=-\int_{0}^\infty \left(-i\sin(\omega \tau)+\tau\e_0 {\bf z}_h(\xi)-\frac{\tau^2}{2}\e_{2n+1}\e_0{\bf z}_h(\xi)^2\right) e^{-p\left({\tau^2{\bf z}_h(\xi)^2+ \frac{\tau^4}{4} {\bf z}_h(\xi)^4}\right)}~\times \\ \ \\ \times \dfrac{E_{\frac{1}{2},\frac{1}{2}}\left(~\cos(\omega \tau) \sqrt{p}~\right)}{\sqrt{p}}~e^{p}dp.
\end{eqnarray*}

Hence, after a wise interchanging on the order of integration of ${\bf K}_\tau(x,t)$, there holds the curious integral representation formula:

\begin{eqnarray*}
	\label{KtauChebyshev}
	{\bf K}_\tau(x,t)=\int_{0}^{\infty} ~\left[ \frac{1}{(2\pi)^{\frac{n}{2}}}  \int_{Q_h} e^{-p\left({\tau^2{\bf z}_h(\xi)^2+ \frac{\tau^4}{4} {\bf z}_h(\xi)^4}\right)} {\mathcal{F}_h{\bf H}_\tau({\bf z}_h(\xi),t;p)}{}e^{-ix\cdot \omega}d\xi\right]~e^{p}dp. \\
\end{eqnarray*}

Here, $\mathcal{F}_h{\bf H}_\tau({\bf z}_h(\xi),t;p)$ stands for the auxiliar function
\begin{eqnarray}
\label{AuxiliarFunction}
\begin{array}{lll}
\mathcal{F}_h{\bf H}_\tau({\bf z}_h(\xi),t;p)&=&\displaystyle -\frac{\tau}{{2\pi}} \int_{-\frac{\pi}{\tau}}^{\frac{\pi}{\tau}}   \left(-i\sin(\omega \tau)+\tau\e_0 {\bf z}_h(\xi)-\frac{\tau^2}{2}\e_{2n+1}\e_0{\bf z}_h(\xi)^2\right)\times \\ \ \\ &\times& \dfrac{E_{\frac{1}{2},\frac{1}{2}}\left(~\cos(\omega \tau) \sqrt{p}~\right)}{\sqrt{p}} ~e^{-i\omega t} d\omega.
\end{array}
\end{eqnarray}

Furthermore, after some straightforward manipulations involving the inversion of the discrete Fourier transform $\mathcal{F}_h$ (cf.~\cite[p.~247]{GuerlebeckSproessig97}) we conclude that the \textit{discrete convolution formula} (\ref{discreteConvolutionformula}) is equivalent to the subordination formula
$$
\Psi(x,t)=
\int_{0}^{\infty} {e^{p\left(-\tau^2\Delta_h+ \frac{\tau^4}{4} \Delta_h^2\right)}}{\bf H}_\tau(D_h,t;p)\left[\Phi_0(x)\right]e^pdp,
$$
involving the \textit{integro-difference operator} ${\bf H}_\tau(D_h,t;p)$ defined through equation (\ref{AuxiliarFunction}), and the semigroup $\left\{e^{p\left(-\tau^2\Delta_h+ \frac{\tau^4}{4} \Delta_h^2\right)}\right\}_{p\geq 0}$ endowed by the Cauchy problem of differential-difference type
\begin{eqnarray}
\label{diffHeatType} \left\{\begin{array}{lll} 
\partial_p \Psi (x,p)= -\tau^2\Delta_h\Psi (x,p)+ \frac{\tau^4}{4} \Delta_h^2 \Psi (x,p) &, (x,p
)\in 
h\BZ^n \times [0,\infty) 
\\ \ \\
\Psi (x,0)=\Phi (x) &, x\in h
\BZ^n
\end{array}\right.. \\ \nonumber
\end{eqnarray}

The description considered above looks tailor-suited for operational purposes, since it combines some well-know facts from Chebyshev polynomials with a fine integral representation involving the \textit{generalized Mittag-Leffler function} (\ref{MittagLeffler}), leading to a quite unexpected link between the Cauchy problem (\ref{DiracCK}) besides the discrete Cauchy-Kovalevskaya extension, and the Cauchy problem (\ref{diffHeatType}) of heat type.
Perhaps an exploitation of this construction may be easily found for Gegenbauer polynomials, or for more general families of ultraspherical polynomials by means of Jacobi expansions (cf.~\cite{Vieira17}). For now we will leave this question open for the interested reader.

Last but not least, 
after completion of this work we learnt from \cite{SteinW2000,MSteinW02,Pierce09} that the incorporation of tools from fractional integral calculus on the discrete setting is almost well-known on the harmonic analysis comunity.
On the Clifford analysis community, this is surely a new research topic that deserves to be developed.

% ------------------------------------------------------------------------
\end{document}